\documentclass[a4paper,12pt]{article}
\usepackage{amsmath,amssymb,amsfonts}
\usepackage{epsfig}
\usepackage{enumerate}
 \usepackage{graphicx}
\begin{document}

\newtheorem{theorem}{Theorem}
\newtheorem{proposition}[theorem]{Proposition}
\newtheorem{lemma}[theorem]{Lemma}
\newtheorem{corollary}[theorem]{Corollary}
\newtheorem{example}[theorem]{Example}
\newtheorem{algo}[theorem]{Algorithm}
\newtheorem{remark}[theorem]{Remark}
\newtheorem{definition}[theorem]{Definition}
\newenvironment{pf}%
{\par\noindent\textbf{Proof:}~}%
{\eop\par\smallskip\par\noindent}

%
%


\newcommand\cJ{\mathcal{J}}
\newcommand\cBC{\mathcal{BC}}
\newcommand\cD{\mathcal{D}}
\newcommand\cE{\mathcal{E}}
\newcommand\bB{\mathcal{BL}}
\newcommand\cB{\mathcal{BL}}
\newcommand\cP{\mathcal{P}}
\newcommand\bF{\mathbf{F}}
\newcommand\cS{\mathcal{S}}
\newcommand\bd{\mathbf{d}}
\newcommand\bc{\mathbf{c}}
\newcommand\bq{\mathbf{q}}
\newcommand\be{\mathbf{e}}
\newcommand\bhate{\hat{ \mathbf{e}}}
\newcommand\bg{\mathbf{g}}

\def\xx{{\bf x}}
\def\aa{{\bf a}}
\def\vv{{\bf v}}
%
\def\cnote#1{\par\noindent\textsf{#1}\par\noindent}
\newcommand\bQ{\mathbf{Q}}
\newcommand\bP{\mathbf{P}}
\newcommand\ba{\mathbf{a}}
\newcommand\bff{\mathbf{f}}
\newcommand\bb{\mathbf{b}}
\newcommand\bh{\mathbf{h}}
\newcommand\bv{\mathbf{v}}
\newcommand\bu{\mathbf{u}}
\newcommand\bw{\mathbf{w}}

\newcommand\cC{\mathcal{C}}
\newcommand\cN{\mathcal{N}}
\newcommand\cR{\mathcal{R}}
\newcommand\ZZ{{\mathbb Z}}
\newcommand{\RR}{\mathbb{R}}
\newcommand{\NN}{\mathbb{N}}
\newcommand{\CC}{\mathbb{C}}
\newcommand{\tN}{{\tilde N}}

\newtheorem{pro}{Proposition}

\newcommand\eg{{\it\thinspace e.g.}}
\def\ie{{\sl \thinspace i.e.},\ }
\newcommand\cpl{\cal PL}
\newcommand{\eop}{\hfill$\Box$}
\newcommand\E{\mathbb{E}}
\def\bla{\mbox{\boldmath $\lambda$}}
\def\bdelta{\mbox{\boldmath $\delta$}}
\def\bmu{\mbox{\boldmath $\nu$}}
\def\bm{\mbox{\boldmath $m$}}
\def\bell{\mbox{\boldmath $\ell$}}
\def\bgamma{\mbox{\boldmath $\gamma$}}
\def\supp{\mbox{\rm supp}\,}

\input epsf.tex
\def\fig#1#2#3{
\vskip #3pc \centerline{{\bf Fig.~#1.} {\it #2}}  }

\def\psone#1#2{\centerline{
    \epsfxsize #2pc  \epsfbox{#1}
    }}

\def\pstwo#1#2#3#4{   \centerline{
\epsfxsize=#2pc \epsfbox{#1}
      \qquad \epsfxsize=#4pc  \epsfbox{#3}
   }}%

\def\psthree#1#2#3#4#5#6{
   \centerline{
    \epsfxsize=#2pc \epsfbox{#1}
      \qquad \epsfxsize=#4pc  \epsfbox{#3}
      \qquad \epsfxsize=#6pc  \epsfbox{#5}
   }}%

\title{\bf Convergence of univariate non-stationary subdivision schemes via asymptotical similarity }
\author{C. Conti\thanks{Costanza Conti - Dipartimento di Ingegneria Industriale,
 Firenze, Italy, {\tt costanza.conti@unifi.it }}, 
 N. Dyn \thanks{Nira Dyn -  School of Mathematical Sciences, Tel-Aviv University,
 Ramat-Aviv, Israel,  {\tt niradyn@post.tau.ac.il}}, C. Manni
\thanks{Carla Manni - Dipartimento di Matematica, Universit\`a di Roma ``Tor Vergata", Italy,
 {\tt manni@axp.mat.uniroma2.it}}, M.-L. Mazure
\thanks{Marie-Laurence Mazure - Laboratoire Jean Kuntzmann, Universit\'e Joseph Fourier, Grenoble, France,
 {\tt Marie-Laurence.Mazure@imag.fr}}}



%
%
\maketitle

\begin{abstract}
\noindent A new equivalence notion between non-stationary subdivision schemes,
termed \emph{asymptotical similarity}, which is weaker than asymptotical equivalence, is introduced and studied. It is known that asymptotical equivalence between a non-stationary subdivision scheme and a convergent stationary scheme guarantees the convergence of the non-stationary
scheme. We show that for non-stationary schemes reproducing constants, the
condition of asymptotical equivalence can be relaxed to asymptotical similarity. This result applies to a wide class of non-stationary schemes of importance in theory and applications. 
\end{abstract}

\bigskip \noindent{{\bf Keywords:}  Non-stationary subdivision schemes,
convergence, reproduction of constants, asymptotical equivalence, asymptotical
similarity 
}


\section{Introduction}
\label{sect: intro}
This short paper studies univariate binary \emph{non-stationary} uniform subdivision
schemes. Such  schemes are efficient iterative methods for genera\-ting smooth
functions via the specification of an initial set of discrete data
$\bff^{[0]}:=\{f^{[0]}_i\in \RR,\ i\in \ZZ\}$,
and a set of refinement rules, 
mapping at each ite\-ration the sequence of values
$\bff^{[k]}:=\{f^{[k]}_i\in\RR,\, i\in \ZZ\}$ attached to the points
of the grid $2^{-k}\ZZ$
into the sequence of values $\bff^{[k+1]}$ attached to the points of
$2^{-(k+1)}\ZZ$.
At each level $k$, the refinement rule $S_{\ba^{[k]}}$, is defined by a
finitely supported \emph{mask} $\ba^{[k]}:=\{a_i^{[k]},\, i\in \ZZ \}$, so that
\begin{equation}\label{def:suboperator}
\bff^{[k+1]}:=S_{\ba^{[k]}}\bff^{[k]}\quad \hbox{with} \quad \left(S_{\ba^{[k]}}\bff^{[k]}\right)_i:=\sum_{j\in \ZZ}a^{[k]}_{i-2j}f_j^{[k]}.
\end{equation}
Each subdivision scheme $\{S_{\ba^{[k]}},\ k\ge 0\}$ we will deal with is
assumed to be  \emph{local}, in the sense that there exists a positive integer
$N$ such that $\supp(\ba^{[k]}):=\{i\in \ZZ,\ |\  a_i^{[k]}\neq0\}\subseteq
[-N,N]$ for all $k\geq 0$.

\medskip \noindent The idea of proving the convergence of a
non-stationary scheme by compa\-rison with a convergent stationary one was
first developed in \cite{DynLevin95}, via the notion of \emph{asymptotical
equivalence} between
non-stationary schemes. Two subdivision schemes $\{S_{\ba^{[k]}},\, k\ge 0\}$
and $\{S_{\ba^{*[k]}},\, k\ge 0\}$ are said to be
asymptotically equivalent when
$$
\sum_{k=0}^\infty  \| S_{\ba^{[k]}}- S_{\ba^{*[k]}}\|<+\infty,
$$
which holds if and only if $\sum_{k=0}^\infty  \| \ba^{[k]}- \ba^{*[k]}\|<
+\infty$.
The main result of the present work is that  for convergence analysis of
non-stationary
schemes reproducing constants, asymptotical equivalence can be replaced by the
weaker notion of \emph{asymptotical similarity}.
We say that two schemes are asymptotically similar when
\begin{equation}\label{eq:UWE}
\lim_{k\rightarrow \infty} \| \ba^{[k]}- \ba^{*[k]}\|=0.
\end{equation}
The class of subdivision schemes to which our result applies is  wide and important from the application point of view. For instance, this class  contains all uniform subdivision schemes generating spaces of exponential polynomials with one exponent equal to zero, and in particular all subdivision schemes for  uniform splines in such spaces \cite{DLL, ContiRomani11}. Besides  their classical interest  in  geo\-metric modelling and approximation
theory, uniform exponential B-splines   
are very useful  in Signal Processing \cite{MU1, MU2} and  in  Isogeometric  Analysis \cite{IgA2, IgA3}. In the latter context, exponential B-splines based subdivision schemes permit to successfully address the difficult evaluation of these splines.

The article is organised as follows. In Sections 2 and 3 the analysis leading to the main result of this paper is
presented. In Section 2  we derive a sufficient condition for the
convergence of
non-stationary schemes reproducing constants, in terms of difference schemes.
This condition replaces the well-known necessary and sufficient condition for
convergence in the stationary case. In Section 3 we introduce the asymptotic similarity relation (\ref{eq:UWE}) and develop some useful consequences for the analysis of non-stationary subdivision schemes.  In particular, we show that, if two subdivision schemes reproduce constants,  and if one of them satisfies the above-mentioned sufficient condition, so does the other. This fact is important for the proof of the convergence of non-stationary schemes reproducing constants by comparison (in the sense of (\ref{eq:UWE})) with convergent stationary ones.  Finally, in Section \ref{sect: EX} we illustrate our result  with non-stationary versions of the de Rham algorithm. \par
\medskip 
Throughout the article the notation $\| \cdot \|$   refers to
the sup-norm, for either operators, functions, or sequences
in $\RR^\ZZ$ and, in particular, we recall that $
 \|S_{\ba^{[k]}}\|\!:=\!\max\!\left(\displaystyle\sum_{i\in \ZZ}\vert a^{[k]}_{2i}\vert,
\sum_{i\in \ZZ}\vert a^{[k]}_{2i+1}\vert\right)$.


\section{A sufficient condition for convergence}
\label{sect: SC}

Let $\{S_{\ba^{[k]}},\ k\ge 0\}$ be a given subdivision scheme, defining successive $\bff^{[k]}$, $k\geq 0$, via (\ref{def:suboperator}). At any level $k\geq 0$, we denote by ${\cpl}(\bff^{[k]})$ the piecewise linear function interpolating  the sequence $\bff^{[k]}$, \ie ${\cpl}(\bff^{[k]})(i2^{-k})=f_i^{[k]}$ for all $i\in \ZZ$. The scheme is said to be \emph{convergent} if, for any bounded $\bff^{[0]}$, the sequence ${\cpl}(\bff^{[k]})$ is uniformly convergent on $\RR$. If so, the limit function is denoted by $S^\infty_{\{\ba^{[k]},\,k\ge 0\}} \bff^{[0]}$.

\smallskip \noindent The subdivision scheme can equivalently be defined by its sequence of symbols,
 the \emph{symbol} of the mask $\ba^{[k]}$ of level $k$ being defined as the Laurent polynomial $a^{[k]}(z):=\sum_{i\in \ZZ}a^{[k]}_{i}z^i$. The scheme $\{S_{\ba^{[k]}},\ k\ge 0\}$ is said to \emph{reproduce constants} if $f^{[0]}_i=1$ for all $i\in\ZZ$ implies  $f^{[k]}_i=1$ for all $i\in\ZZ$ and all $k\geq 0$, which holds if and only if
$$
\sum_{i\in \ZZ}a^{[k]}_{2i}=\sum_{i\in \ZZ}a^{[k]}_{2i+1}=1\quad \hbox{for all }k\ge 0,
$$
or if and only if the symbols satisfy
\begin{equation}\label{eq:repcsymb}
a^{[k]}(-1)=0 \hbox{ and } a^{[k]}(1)=2\quad \hbox{for all } k\ge 0.
\end{equation}
If (\ref{eq:repcsymb}) holds, each symbol can be written as $a^{[k]}(z)=(1+z)q^{[k]}(z)$, where $q^{[k]}(z):=\sum_{i\in \ZZ}q^{[k]}_{i}z^i$  satisfies $q^{[k]}(1)=1$, and we have
\begin{equation}\label{eq:q}
q_{i}^{[k]}=\sum_{j\leq i}(-1)^{i-j}a_{j}^{[k]}, \quad a_i^{[k]}=q_i^{[k]}+ q_{i-1}^{[k]}, \quad i\in\ZZ, \ k\geq 0.
\end{equation}
From the rightmost relation in (\ref{eq:q}) it is easily seen that the scheme  $\{S_{\bq^{[k]}},k\ge 0\}$ permits the computation of all backward differences
 $\Delta f_i^{[k]}:=f_i^{[k]}-f_{i-1}^{[k]}$, namely  
$$
 \Delta \bff^{[k+1]}= S_{\bq^{[k]}} \Delta \bff^{[k]},   \quad\hbox{with }\Delta \bff^{[k]}:=\{\Delta f^{[k]}_i,i\in \ZZ\}. 
 $$
The non-stationary subdivision scheme $\{S_{\bq^{[k]}},k\ge 0\}$ is called the difference scheme of $\{S_{\ba^{[k]}},k\ge 0\}$.

\smallskip 
The scheme $\{S_{\ba^{[k]}},k\ge 0\}$ is \emph{stationary}  when its masks $\ba^{[k]}$ do not depend on the level $k$, \ie $\ba^{[k]}=\ba$ for all $k\geq 0$. In that case we will use the simplified notation $\{S_\ba\}$. 

As is well known, reproduction of constants is necessary for convergence of stationary subdivision schemes. Let us also recall the following other major fact of the stationary case (see \eg \cite{Dyn2002}).

\begin{theorem}\label{the: ContractStat}
 Let $\{S_\ba\}$ be a stationary subdivision scheme reproducing constants, with difference scheme $\{S_\bq\}$. Then
 the scheme  $\{S_\ba\}$ converges if and only if there exists a positive integer $ n$ such that $\mu:=\left\| (S_{\bq})^n\right\| < 1$.
\end{theorem}
 \noindent
A similar necessary and sufficient condition for the convergence of non-stationary subdivision schemes is not known. Nevertheless, a non-stationary version of  
the sufficient condition is given in Theorem \ref{SCC} below.
\begin{definition}
We say that a subdivision scheme $\{S_{\ba^{[k]}}, k\ge 0\}$, assumed to reproduce constants, satisfies \emph{Condition A}, when its difference scheme $\{S_{\bq^{[k]}}, k\ge 0\}$ fulfills  the following requirement: \begin{equation}\label{eq:scc}
	\begin{split}
&\hbox{there exist two integers } K\geq 0,\  n>0,  \hbox{ such that}\\
&\quad\mu:=\displaystyle\sup_{k\geq K}\left\| S_{\bq^{[k+n-1]} }\ldots S_{\bq^{[k+1]} }S_{\bq^{[k]} }\right\| < 1.
	\end{split}
\end{equation}
\end{definition}

\noindent
Let us recall that a scheme $\{S_{\ba^{[k]}},k\ge 0\}$ is said to be \emph{bounded}, if $\sup_{k\geq 0}\|S_{\ba^{[k]}} \|<+\infty$, or, equivalently, due to locality, if $\sup_{k\geq0}\|\ba^{[k]} \|<+\infty$.

\begin{theorem}
\label{SCC}
Let $\{S_{\ba^{[k]}},k\ge 0\}$ be a bounded subdivision scheme reprodu\-cing constants and satisfying Condition A. Then,  $\{S_{\ba^{[k]}},k\ge 0\}$ converges.  Moreover, there exists a positive number $C$, such that, for any initial $\bff^{[0]}$,
\begin{equation}\label{eq:uconv}
\|S^\infty_{\{\ba^{[k]},\,k\ge 0\}} \bff^{[0]}-{\cpl}\left(\bff^{[k]}\right)\|\le C\ \widehat\mu^k \|\Delta \bff^{[0]}\|, \quad k\ge 0, \quad \hbox{with } \widehat\mu:=\mu^{  \frac{1}{n}},
\end{equation}
where $\mu$ and $n$ are provided by (\ref{eq:scc}), and where $\{\bff^{[k]},\ k\ge 0\}$ are the sequences generated by the subdivision scheme.
\end{theorem}

\smallskip \noindent Before proving the theorem we prove two lemmas. Below, as well as whenever we refer to a specific mask, we only indicate the non-zero elements.

 \begin{lemma}\label{lemma:preliminary}
Let $\{S_{\ba^{[k]}}, k\ge 0\}$ be a bounded subdivision scheme which reproduces constants, its locality being prescribed by the positive integer $N$. Let $\bh:=\left\{\frac12,1,\frac12\right\}$ be the mask of the stationary linear B-spline subdivision scheme. The symbols of the masks $\{\bd^{[k]}\!:=\ba^{[k]}\!-\!\bh,k\ge 0\}$ can  be written as
\begin{equation}\label{eq:d}
d^{[k]}(z)=(1-z^2)e^{[k]}(z),
\end{equation}
where, for each $k\geq 0$, the mask $\be^{[k]}$ satisfies
\begin{equation}\label{eq:e}
e_{i}^{[k]}:=\sum_{j\geq 0}d_{i-2j}^{[k]}, \quad \hbox{for all\ }i\in\ZZ, \quad \supp \be^{[k]}\subset [-N, N-2].
\end{equation}
\end{lemma}
\pf
The factorization (\ref{eq:d}) is valid for the difference of any two subdivision schemes reproducing constants since their symbols take the same value at $-1$ and $1$, see (\ref{eq:repcsymb}). The rest of the claim readily follows from \eqref{eq:d}.
\eop

\begin{lemma}\label{lemma:W}
Under the assumptions of Theorem \ref{SCC}
there exists a positive constant $C_1$ such that
\begin{equation}\label{cond-I}
    \Vert\Delta \bff^{[k]}\Vert\le C_1\ \widehat \mu^k \Vert\Delta \bff^{[0]}\Vert, \quad k\geq 0. 
    \end{equation}
\end{lemma}

\pf
Select any integers $p,r$, with $p\geq 0$ and $0\leq r\leq n-1$, where $n$ is given by \eqref{eq:scc}. Repeated application of (\ref{eq:scc}) yields:
\begin{equation}\label{cond-Ibis}
	\begin{split}
    \Vert\Delta \bff^{[K+pn+r]}\Vert&\le  \mu^p \Vert\Delta \bff^{[K+r]}\Vert, \\
                                                          &\le \widehat\mu^{K+pn+r} \frac{\Vert S_{\bq^{[K+r-1]}}\ldots S_{\bq^{[1]}}S_{\bq^{[0]}}\Vert}{\widehat \mu^{K+r} }\Vert\Delta \bff^{[0]}\Vert.
	\end{split}
\end{equation}
From (\ref{cond-Ibis}) and from the fact that $\widehat \mu<1$ it can easily be derived that (\ref{cond-I}) holds with
$$
C_1:=\frac{1}{\widehat\mu^{K+n-1}}\max_{0\leq k\leq K+n-1}\Vert S_{\bq^{[k-1]} }\ldots S_{\bq^{[1]} }S_{\bq^{[0]} }\Vert.
$$\eop

\noindent
{\bf Proof of Theorem \ref{SCC}:}
By standard arguments
it is sufficient to show that the sequence $F^{[k]}:={\cpl}(\bff^{[k]})$, $k\ge 0$, of piecewise linear interpolants satisfies
\begin{equation}\label{eq:F-F}
\|F^{[k+1]}-F^{[k]}\|\le \Gamma\  \widehat \mu^k\Vert\Delta \bff^{[0]}\Vert\,, \quad k\geq 0,
\end{equation}
for some positive constant $\Gamma$. The constant $C$ in (\ref{eq:uconv}) can then be chosen as $C:=\Gamma/(1- \widehat \mu)$. With the help of the hat function $$H(x)=\left\{
                                  \begin{array}{ll}
                                    1-|x|, & x\in (-1,1), \\
                                    0, & \hbox{otherwise},
                                  \end{array}
                                \right.
$$
we can write $F^{[k+1]}$ and $F^{[k]}$ respectively as
$$
F^{[k+1]}(x)=\sum_{i\in \ZZ}\left(S_{\ba^{[k]}}\bff^{[k]}\right)_iH(2^{k+1}x-i)\,,
$$
and
$$
F^{[k]}(x)=\sum_{i\in \ZZ}\bff^{[k]}_i H(2^kx-i)=\sum_{i\in \ZZ}\left(S_{\bh}\bff^{[k]}\right)_iH(2^{k+1}x-i)\,,
$$
where $S_{\bh}$ is the subdivision scheme for linear B-splines recalled in Lemma \ref{lemma:preliminary}. Hence, by the definition of $\bd^{[k]}$ in Lemma \ref{lemma:preliminary}, we obtain
\begin{equation}\label{eq:Fk+1-Fk}
F^{[k+1]}(x)-F^{[k]}(x)=\sum_{i\in \ZZ}g^{[k+1]}_iH(2^{k+1}x-i)  \hbox{ with }\bg^{[k+1]}:=S_{\bd^{[k]}}\bff^{[k]}.
\end{equation}
The left relations in (\ref{eq:e}) can be written as
$d_i^{[k]}=e_i^{[k]}-e_{i-2}^{[k]}$ for all $i\in\ZZ$,
implying that
\begin{equation}\label{eq:g}
g^{[k+1]}_i=\sum_{j\in \ZZ}e^{[k]}_{i-2j}\left(\Delta\bff^{[k]}\right)_j,  \quad i\in\ZZ.
\end{equation}
Now, Lemma \ref{lemma:preliminary} and the boundedness assumption ensure that
\begin{equation}\label{eq:C2}
\Vert \be^{[k]}\Vert \leq  C_2:=N(\sup_{j\geq 0}\Vert \ba^{[j]}\Vert+1)<+\infty, \quad k\geq 0.
\end{equation}
Gathering (\ref{eq:C2}), (\ref{eq:g}), (\ref{eq:Fk+1-Fk}), (\ref{cond-I}) leads to (\ref{eq:F-F}), with $\Gamma:=NC_1C_2$.
 \eop

\bigskip \noindent
As in the stationary case, it can be proved that the limit function in Theorem \ref{SCC} is H\"{o}lder continuous with exponent $|Log_2 \widehat\mu|$.

\begin{remark}
{\rm Different proofs of the fact that Condition A is sufficient for convergence  already exist in the wider context of non-regular (i.e, non-uniform, non-stationary) schemes, using non-regular grids,  either nested \cite{MaximMazure2004} or non-nested \cite{MLM2005, MLM2006}. Nevertheless, we did consider it useful to give a simplified proof  in the context of uniform schemes and  regular grids. Indeed, in that case the proof is made significantly more accessible by the use of the corresponding classical tools.   }
\end{remark}

\section{Asymptotically similar schemes}
\label{sect: AS}

\begin{definition}
We say that  two subdivision schemes $\{S_{\ba^{[k]}},\, k\ge 0\}$ and $\{S_{\ba^{*[k]}}$, $k\ge 0\}$ are
\emph{asymptotically similar} if they satisfy
\begin{equation}\label{eq:equiv}
\lim_{k\rightarrow\infty}\|\ba^{[k]}-\ba^{*[k]}\|=0.
\end{equation}
\end{definition}
Clearly, asymptotical similarity is an equivalence relation between subdivision schemes, which is weaker than asymptotical equivalence, see \cite{DynLevin95}.  By the locality of the two schemes,  proving  their asymptotical similarity simply consists in checking that
$$
\lim_{k\rightarrow\infty}(a^{[k]}_i-a^{*[k]}_i)=0\quad\hbox{for }-N\leq i\leq N\,,
$$
where $[-N,N]$ contains the support of the masks $\ba^{[k]},\ \ba^{*[k]}$ for $k\ge 0$.
 Note that (\ref{eq:equiv}) can be replaced by $\lim_{k\rightarrow \infty} \| S_{\ba^{[k]}}-S_{\ba^{*[k]}}\|=0$ as well. If two  subdivision schemes are asymptotically similar and if one of them is bounded, so is the other. 
 
\smallskip \noindent Depending on the properties of the schemes, asymptotical similarity can be expressed in different ways:
\begin{proposition}\label{prop: equiv}
Given two subdivision schemes $\{S_{\ba^{[k]}},\, k\ge 0\}$ and $\{S_{\ba^{*[k]}},\, k\ge 0\}$ which both reproduce constants, the following properties are equivalent:
\begin{enumerate}[\rm(\roman{enumi})]
\item $\{S_{\ba^{[k]}},\,k\ge 0\}$ and $\{S_{\ba^{*[k]}},\,k\ge 0\}$ are asymptotically similar;
\item the difference schemes $\{S_{\bq^{[k]}},\,k\ge 0\}$ and $\{S_{\bq^{*[k]}},\,k\ge 0\}$ are asymptotically similar.
\end{enumerate}
If, in addition, one of the two subdivision schemes $\{S_{\ba^{*[k]}},\, k\ge 0\}$ or $\{S_{\ba^{[k]}},\, k\ge 0\}$ is bounded, then {\rm(i)} is also equivalent to
\begin{enumerate}[\rm(\roman{enumi})]
\setcounter{enumi}{2}
\item for any fixed $p\geq 0$, $\lim_{k\rightarrow\infty}\left\|S_{\bq^{[k+p]} }\ldots S_{\bq^{[k]} }-S_{\bq^{*[k+p]} }\ldots S_{\bq^{*[k]} }\right\|=0$.
\end{enumerate}
\end{proposition}

\pf Without loss of generality we can assume that the locality of the two schemes is determined by the same positive  integer $N$. Then, by application of (\ref{eq:q}) we can derive that
$\supp \bq^{[k]},\ \ \supp \bq^{*[k]}\subset [-N, N-1]$, and that
$$
\frac{1}{2}\| \ba^{[k]}-\ba^{*[k]}\|\leq \| \bq^{[k]}-\bq^{*[k]}\|\leq 2N\| \ba^{[k]}-\ba^{*[k]}\|\,.
$$
The equivalence between (i) and (ii) follows. Clearly, (ii) is implied by (iii). As for the implication $\rm(ii)\Rightarrow(iii)$, it follows by induction from  the equality
$${
\begin{split}
S_{\bq^{[k+p+1]} }&S_{\bq^{[k+p]} }\ldots S_{\bq^{[k]} }-S_{\bq^{*[k+p+1]} }S_{\bq^{*[k+p]} }\ldots S_{\bq^{*[k]}}\\
&=
\left(S_{ \bq^{[k+p+1]} }-S_{ \bq^{*[k+p+1]} }\right)S_{\bq^{[k+p]} }\ldots S_{\bq^{[k]} }\\
&\qquad+ S_{ \bq^{*[k+p+1]} }\left(S_{\bq^{[k+p]} }\ldots S_{\bq^{[k]} }-S_{\bq^{*[k+p]} }\ldots S_{\bq^{*[k]} }\right), \quad k\geq 0\,,
\end{split}
}
$$
and from the boundedness of the two schemes. \eop

\begin{proposition}\label{prop: condA}
Let $\{S_{\ba^{*[k]}},k\ge 0\}$ be a bounded subdivision scheme reproducing constants and satisfying Condition A. Then, any subdivision scheme $\{S_{\ba^{[k]}},k\ge 0\}$ which reproduces constants and is asymptotically similar to $\{S_{\ba^{*[k]}}, k\ge 0\}$, also satisfies Condition A.
\end{proposition}

\pf
We know the existence of two integers $K^*, n$ such that
$$
\mu^*:=\displaystyle\sup_{k\geq K^*}\left\| S_{\bq^{*[k+n-1]} }\ldots S_{\bq^{*[k]} }\right\| < 1.
$$
Select any $\mu\in(\mu^*, 1)$ and choose $\varepsilon>0$ such that $\mu^*+\varepsilon<\mu$. The two schemes being asymptotically similar,  and $\{S_{\ba^{*[k]}},k\ge 0\}$ being bounded, we know that (iii) of Proposition \ref{prop: equiv} holds. We can thus find $\widetilde K\geq 0$, such that
$$
\left\|S_{\bq^{[k+n]} }\ldots S_{\bq^{[k]} }-S_{\bq^{*[k+n]} }\ldots S_{\bq^{*[k]} }\right\|\leq \varepsilon \hbox{ for all }k\geq \widetilde K.
$$
Clearly, we have
\begin{equation}\label{eq:mu}
\left\|S_{\bq^{[k+n]} }\ldots S_{\bq^{[k]} }\right\| \leq \mu <1, \quad \hbox{for each }k\geq K:=\max(K^*,  \widetilde K).
\end{equation}
The claim is proved.
\eop

\begin{remark}  {\em We would like to draw the reader's attention to the fact that we have not proved that, when two bounded non-stationary subdivision schemes reproducing constants are asymptotically similar, convergence  of one of them implies  convergence of the other.  Convergence of the second scheme is obtained only when convergence of the first one results from Condition A. This follows from Proposition \ref{prop: condA} and Theorem \ref{SCC}. This is actually sufficient to prove Theorem \ref{teo:main}Ê below, which is the main application of all previous results. 
} 
\end{remark}

\begin{theorem}\label{teo:main}
Let $\{S_{\ba^*}\}$ be a convergent stationary subdivision scheme with $\mu^*:=\left\|(S_{\bq^*})^n\right\|<1$. Let $\{S_{\ba^{[k]}},\ k\ge 0\}$ be a non-stationary subdivision scheme reproducing constants which  is asymptotically similar to $\{S_{\ba^*}\}$.
Then, the  scheme $\{S_{\ba^{[k]}},\ k\ge 0\}$ is convergent and for any $\eta \in ({\mu^*}^{\frac1 n},1)$ there exists a positive constant $C$ such that, for any initial bounded $ \bff^{[0]}$,
$$
\|S^\infty_{\{\ba^{[k]},\ k\ge 0\}} \bff^{[0]}-{\cpl}\left(\bff^{[k+1]}\right)\|\le C\ \eta^k  \|\Delta \bff^{[0]}\|, \quad k\ge 0.
$$
\end{theorem}

\pf
The existence of a positive integer  $n$ with $\mu^*:=\left\|(S_{\bq^*})^n\right\|<1$ is due to  the stationary scheme $\{S_{\ba^*}\}$ being convergent, see Theorem \ref{the: ContractStat}. In other words, $\{S_{\ba^*}\}$ satisfies Condition A. We also know that $\{S_{\ba^*}\}$ reproduces constants.  Accordingly,  by application of Proposition \ref{prop: condA}, we can say that $\{S_{\ba^{[k]}}, k\geq 0\}$  satisfies Condition A too. Furthermore, we know that we can apply Theorem \ref{SCC} using any $\mu\in (\mu^*, 1)$ (see (\ref{eq:mu})). \eop

\section{Illustrations}
\label{sect: EX}

In order to illustrate the usefulness of asymptotic similarity, in particular via Theorem \ref{teo:main},  we consider a non-stationary version of the de Rham algorithm.
\begin{figure}
\begin{center}
\includegraphics[width=7cm]{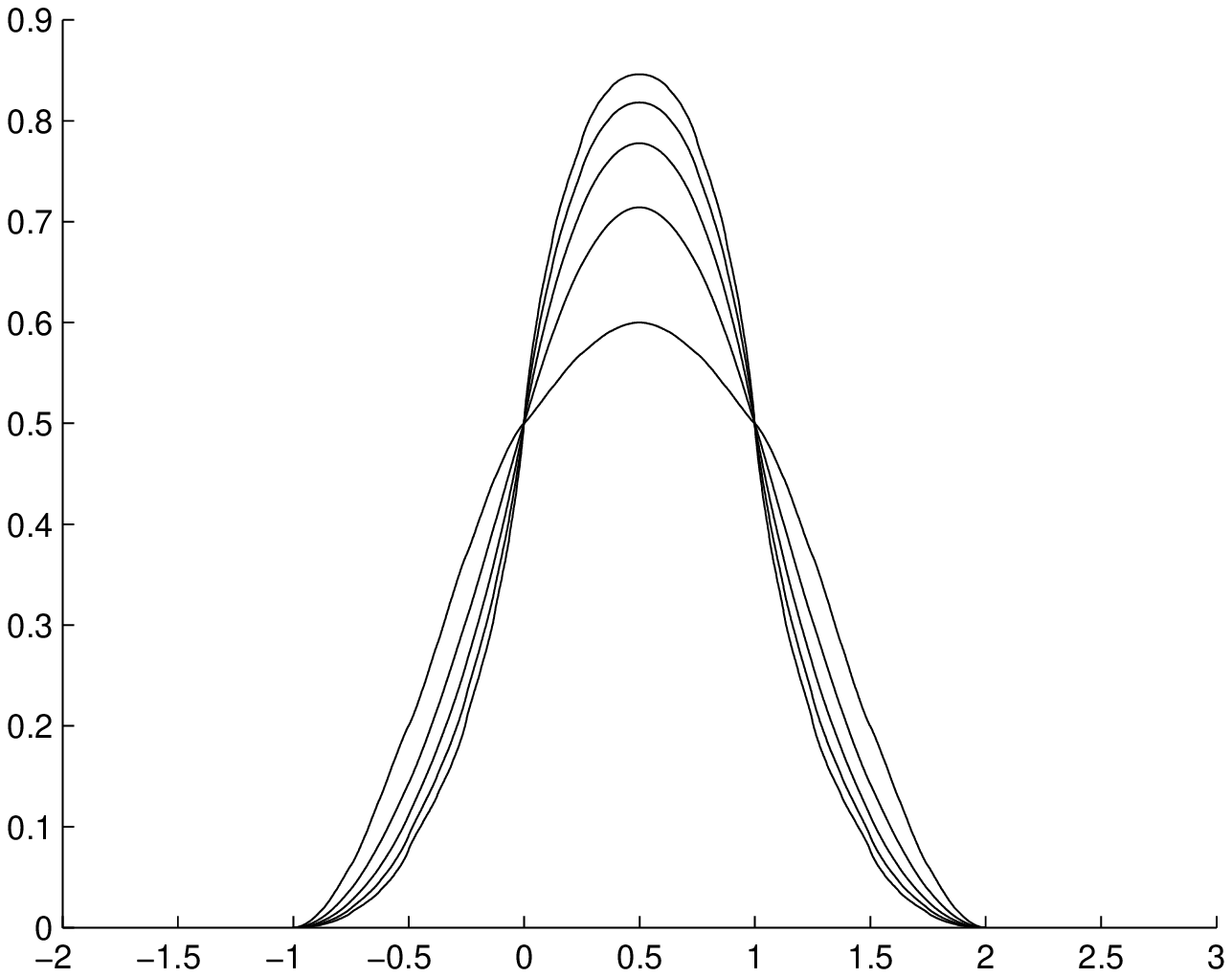}\hskip -.5cm\includegraphics[width=7cm]{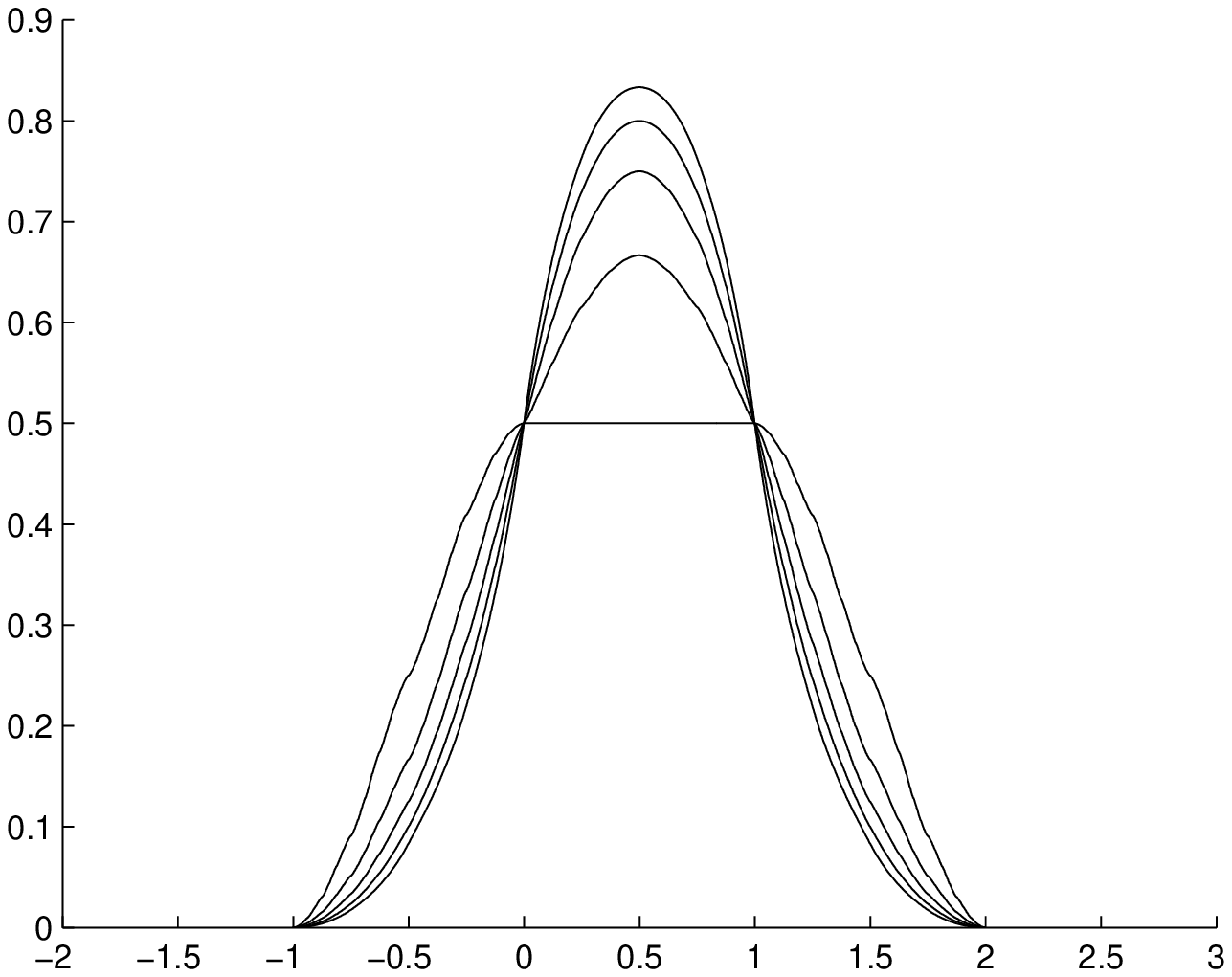}
\caption{\small Limit functions obtained via (\ref{eq:rham}) and (\ref{eq:rhambis}) starting from $\bff^{[1]}=\bdelta$, with $\gamma=2$ (left) and $\gamma=1.5$ (right), and with, everywhere $\varepsilon_k=\frac{\alpha}{k}$, $k\geq 1$. For both pictures  the five displayed functions correspond to $\alpha=$ 2.5; 1.5;  0.5; -0.5;  -1.5  (from top to bottom).}
\end{center}
\end{figure}
\begin{figure}
\begin{center}
\vskip-.5cm
\includegraphics[width=4.8cm]{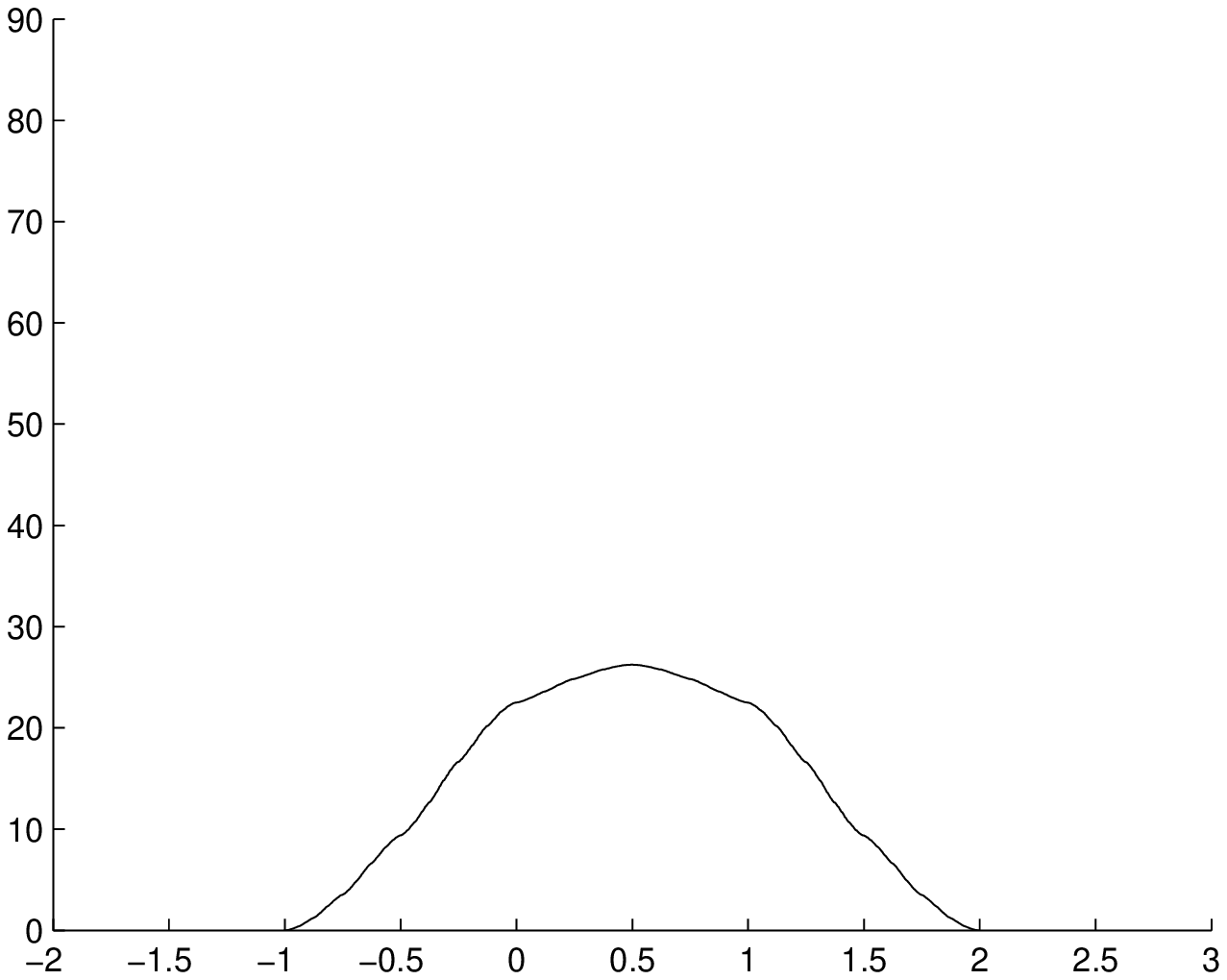}
\hskip-.5cm\includegraphics[width=4.8cm]{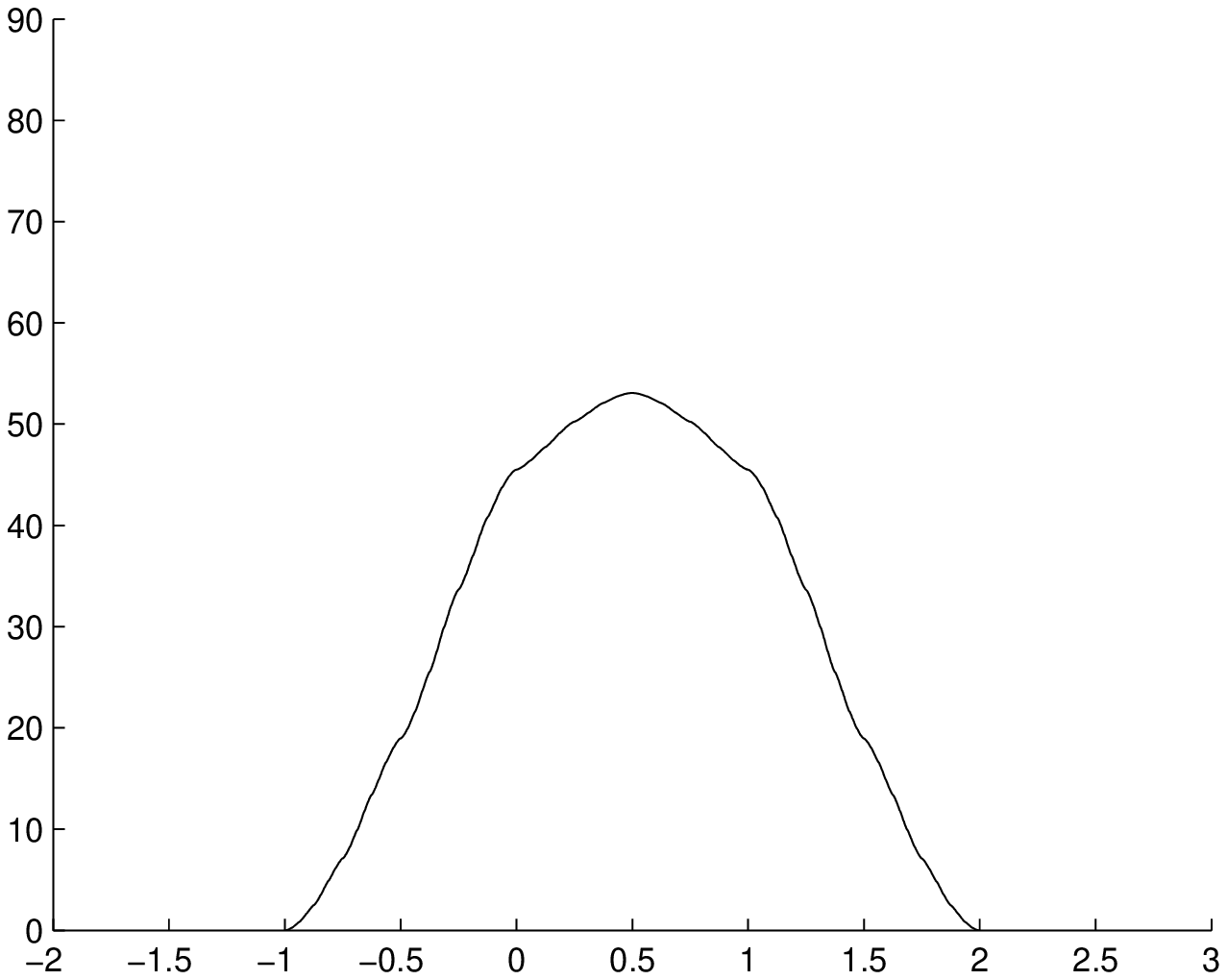}\hskip-.5cm\includegraphics[width=4.8cm]{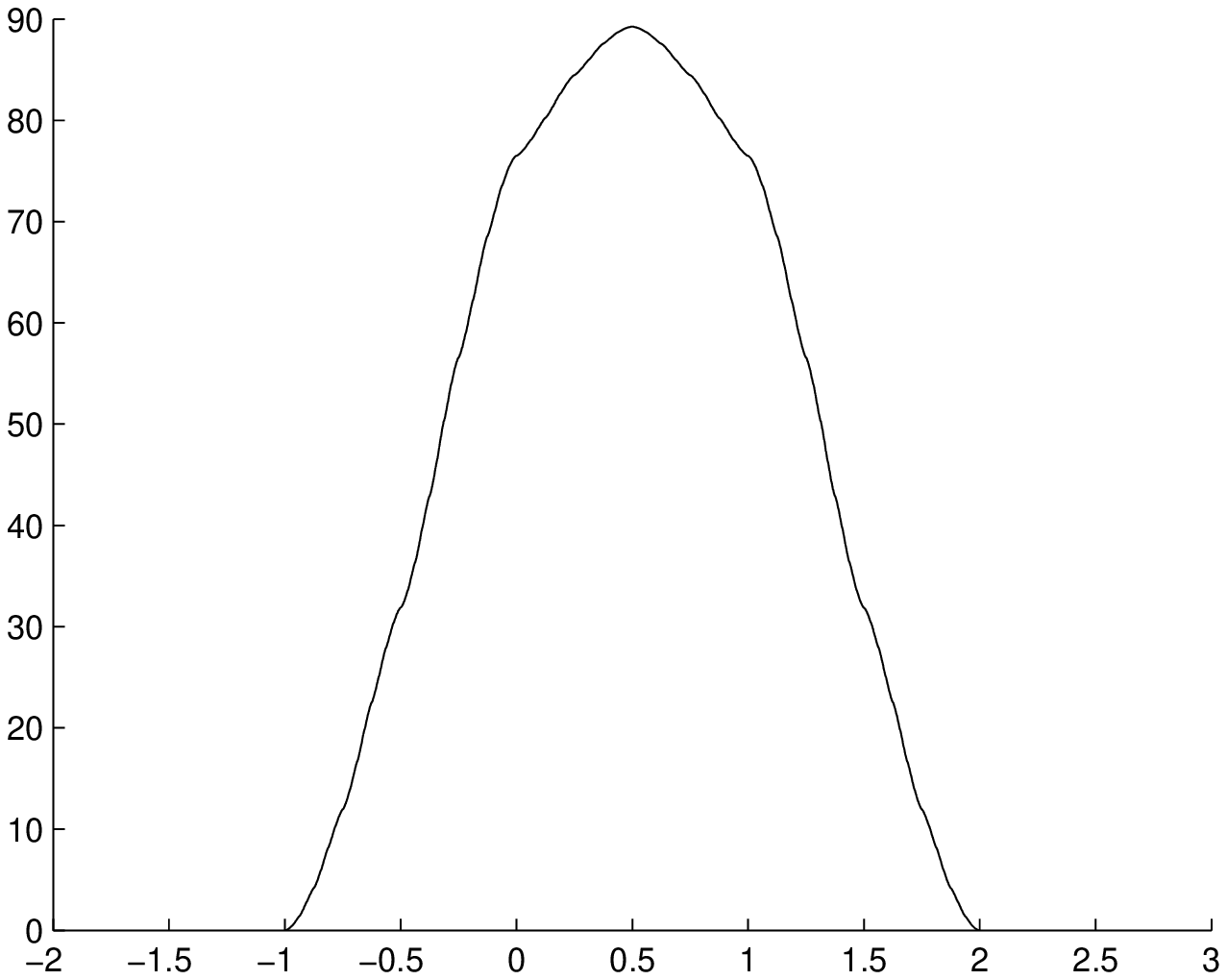}
\caption{\small From left to right: 8; 12; 16 iterations of  (\ref{eq:chaikinter}) starting from $\bff^{[1]}=\bdelta$. }
\end{center}
\end{figure}

At each level $k$, exactly two consecutive points of the next level are located on each segment of the polygonal line of level $k$, so that they divide the segment with ratios $1:\gamma_k:1$, where $\gamma_k$, $k\geq 0$,  is a given sequence of positive numbers. We obtain a non-stationary subdivision scheme $\{S_{\ba^{[k]}},\ k\ge 0\}$ with masks
\begin{equation}\label{eq:rham}
\ba^{[k]}=\left\{\underbrace{\frac{1}{2+\gamma_k}}_{a_0^{[k]}},\  \frac{1+\gamma_k}{2+\gamma_k}, \ \frac{1+\gamma_k}{2+\gamma_k}, \ \frac{1}{2+\gamma_k}\right\}, \quad k\geq 0.
\end{equation}
We additionally assume the existence of a positive number $\gamma$ such that
\begin{equation}\label{eq:rhambis}
\gamma_k=\gamma+\varepsilon_k\hbox{ for all }k\geq 0, \quad \hbox{with } \lim_{k\rightarrow \infty} \varepsilon_k=0.
\end{equation}
\noindent

The non-stationary subdivision scheme defined by (\ref{eq:rham}) and (\ref{eq:rhambis}) is asymptotically similar to the classical stationary de Rham scheme $\{S_{\ba^*}\}$ which is obtained when all
$\gamma_k$'s are equal to $\gamma$ \cite{rham} (see also \cite{ContiRomani10} and \cite{ContiRomani13}). Indeed, all masks have the same support and
$$
a_i^*-a_i^{[k]}=\pm\frac{\varepsilon_k}{(2+\gamma_k)(2+\gamma)} \hbox{ for all }i\in\supp(\ba^{[k]}) \hbox{ and for all }k\geq 0.
$$
Since all assumptions of Theorem \ref{teo:main} are satisfied, $\{S_{\ba^{[k]}},\ k\ge 0\}$ converges when  $\{S_{\ba^*}\}$ converges, that is, for all positive $\gamma$. 
We illustrate this in Figure 1, where, for  $\gamma=2$, and $\gamma=1.5$, limit functions corresponding to various sequences $\varepsilon_k$, $k\geq 1$, are shown, starting from the initial sequence $\bff^{[1]}:=\bdelta=\{\delta_{i,0},\ i\in \ZZ\}$. For $\gamma=2$, $\{S_{\ba^*}\}$ is simply the Chaikin algorithm with mask $\ba=\{\frac14,\ \frac 34,\ \frac 34,\ \frac 14\}$. In either illustration, the non-stationary subdivision scheme $\{S_{\ba^{[k]}},\ k\ge 1\}$ is not asymptotically equivalent to the corresponding de Rham scheme.

\medskip
For the family of masks $\{\ba^{[k]},\ k\ge 1\}$ with
\begin{equation}\label{eq:chaikinter}
\ba^{[k]}=\left\{\underbrace{\frac 14+\frac1k}_{a_0^{[k]}},\ \frac 34 +\frac1k,\ \frac 34+\frac 1k,\ \frac14+\frac1k\right\},\quad k\ge 1, 
\end{equation}
Figure 2 shows the results after $8,\ 12,\ 16$ iterations in the left, in the center and in the right, respectively. It clearly shows that the corresponding non-stationary scheme is not convergent.
Still, it is asymptotically similar  to the Chaikin scheme as in the scheme in the left side of Figure 1. This is not in contradiction with Theorem \ref{teo:main} since reproduction of  constants is not satisfied. Indeed, $\sum_{i\in \ZZ}a^{[k]}_{2i}=\sum_{i\in \ZZ}a^{[k]}_{2i+1}=1+\frac 2k\neq 1$ for all $k\ge 1$. This enhances the importance of  all assumptions for the validity of Theorem \ref{teo:main}.

\section{Conclusion}
Non-stationary subdivision schemes are not as simple to handle as their stationary counterparts. Analyzing them by comparison with a simpler scheme is quite a natural idea. Up to now, the main tool for such a comparison was the asymptotical equivalence, as developed in \cite{DynLevin95}, see also \cite{DynLevinJoon2007}.  Still, relevant examples show that this is sometimes a too demanding requirement. This motivated the present note, in which we have replaced asymptotical equi\-valence by asymptotical similarity, a simpler and weaker equivalence relation  between non-stationary schemes. Provided that it reproduces constants, a non-stationary scheme which is asymptotically similar to a convergent stationary one is convergent. The proof relies on a sufficient condition for convergence invol\-ving differences schemes.

To enhance the interest of asymptotic similarity, we would like to mention that this notion can be adapted to the non-regular framework where  it yields interesting results, see \cite{MLM2013}. 


\end{document}